\documentclass{amsart}
\title{Weakly infinite dimensional subsets of $\reals^{\naturals}$}
\author{Liljana Babinkostova and Marion Scheepers}

\usepackage{amsfonts}
\newcommand{\sone}{{\sf S}_1}

\newcommand{\sfin}{{\sf S}_{fin}}

\newcommand{\Sc}{{\sf S}_c}

\newcommand{\op}{\mathcal{O}}
\newcommand{\om}{\Omega}
\newcommand{\ga}{\Gamma}

\newcommand{\opk}{\mathcal{O}_{\sf k}}
\newcommand{\opkfd}{\mathcal{O}_{\sf kfd}}

\newcommand{\opcfd}{\mathcal{O}_{\sf cfd}}
\newcommand{\opfd}{\mathcal{O}_{\sf fd}}

\newcommand{\naturals}{{\mathbb N}}
\newcommand{\reals}{{\mathbb R}}

\newcommand{\irrationals}{{\mathbb P}}

\newtheorem{theorem}{{\bf Theorem }}
\newtheorem{proposition}[theorem]{{\bf Proposition }}
\newtheorem{lemma}[theorem]{{\bf Lemma }}
\newtheorem{corollary}[theorem]{{\bf Corollary}}	
\newtheorem{problem}{{\bf Problem }}

\newtheorem{conjecture}{{\bf Conjecture}}
\newtheorem{example}{{\bf Example}}
\begin{document}
\maketitle
\begin{abstract} The Continuum Hypothesis implies an Erd\"os-Sierpi\'nski like duality between the ideal of first category subsets of $\reals^{\naturals}$, and the ideal of countable dimensional subsets of $\reals^{\naturals}$. The algebraic sum of a Hurewicz subset - a dimension theoretic analogue of Sierpinski sets and Lusin sets - of $\reals^{\naturals}$ with any compactly 
 countable dimensional subset of $\reals^{\naturals}$ has first category. 
\end{abstract}

\section{Introduction.}

The main purpose of this paper is threefold: We point out strong analogies between the well-known class of countable dimensional subsets of the Hilbert cube, and the classes of Lebesgue measure zero subsets of the real line and the class of first category subsets of complete separable metric spaces. A classical theorem of Erd\"os and Sierpi\'nski gives, under the Continuum Hypothesis, some explanation of this analogy. We show also that the dimension theoretic analogue of Lusin sets and Sierpinski sets has some of the deeper properties shared by Lusin sets and Sierpi\'nski sets. And we give some information about a class of weakly infinite dimensional spaces which is emerging as an important class in dimension theory. In particular we show, using the Continuum Hypothesis, that the most restrictive of these classes contains metric spaces which are not countable dimensional.

$\reals^{\naturals}$ denotes the Tychonoff product of countably many copies of the real line, $\reals$. The symbol $\lbrack 0,\, 1\rbrack$ denotes the unit interval and the subspace $[0,1]^{\naturals}$ of $\reals^{\naturals}$ denotes the Hilbert cube. The symbol $\irrationals$ denotes the set of irrational numbers, viewed as a subspace of the real line.

\section{The countable dimensional subsets of $\reals^{\naturals}$}

Hurewicz \cite{hurewicz28} defined a subset of a topological space to be \emph{countable dimensional} if it is a union of countably many zero-dimensional subsets of the space, and proved that the Hilbert cube, thus also $\reals^{\naturals}$ is not countable dimensional. Thus, the collection {\sf CD} of countable dimensional subsets of $\reals^{\naturals}$ is a $\sigma$-ideal. Since each metrizable space of cardinality less than $2^{\aleph_0}$ is zerodimensional, the cofinality of {\sf CD}, denoted {\sf cof}({\sf CD}), is $2^{\aleph_0}$. Smirnov and independently Nagami also showed that every separable metric space is a union of $\aleph_1$ zero-dimensional subsets, and thus the covering number of {\sf CD}, denoted {\sf cov}({\sf CD}), is $\aleph_1$. This is unlike the $\sigma$-ideal $\mathcal{M}$ of first category subsets of $\reals$, and the $\sigma$-ideal $\mathcal{N}$ of Lebesgue measure zero subsets of $\reals$, where the cofinality and the covering number are not decided by ZFC. 

If one assumes the Continuum Hypothesis (CH), these differences vanish and one can give an explanation for the unusually strong analogies between the theory of the countable dimensional subspaces of $\reals^{\naturals}$ and the theory of the first category subsets of $\reals^{\naturals}$: These two collections are dual in the sense of Erd\"os and Sierpi\'nski. Chapter 19 of \cite{Ox} gives a nice exposition of the following classical theorem:
\begin{theorem}[Erd\"os-Sierpi\'nski Duality Theorem]\label{erdsierp} Let $X$ be a set of cardinality $\aleph_1$. Let $\mathcal{I}$ and $\mathcal{J}$ be $\sigma$-complete ideals on $X$ such that:
\begin{enumerate}
  \item{$X=A\bigcup B$ where $A$ and $B$ a disjoint sets with $A\in\mathcal{I}$ and $B\in\mathcal{J}$.} 
  \item{${\sf cof}(\mathcal{I}) = \aleph_1 = {\sf cof}(\mathcal{J})$.} 
  \item{For each $I\in\mathcal{I}$ there is an $S\subset X\setminus I$ such that $\vert S\vert=\aleph_1$ and $S\in\mathcal{I}$.} 
  \item{For each $J\in\mathcal{J}$ there is an $S\subset X\setminus J$ such that $\vert S\vert=\aleph_1$ and $S\in\mathcal{J}$.} 
\end{enumerate}
Then there is a bijective function $f:X\longrightarrow X$ which is its own inverse, such that for each set $E\subseteq X$ we have
\[
  E\in\mathcal{I} \Leftrightarrow f\lbrack E\rbrack \in\mathcal{J}.
\]
\end{theorem}

A classical theorem of Tumarkin \cite{tumarkin} is an important tool in proving the duality between countable dimensional sets and first category sets of $\reals^{\naturals}$ under CH:
\begin{theorem}[Tumarkin]\label{tumarkin} In a separable metric space each $n$-dimensional set is contained in an $n$-dimensional ${\sf G}_{\delta}$-set. 
\end{theorem}

\begin{corollary}\label{dimandcat} Assume the Continuum Hypothesis. There is a bijection $f:\reals^{\naturals}\longrightarrow\reals^{\naturals}$ such that $f = f^{-1}$ and for each set $E\subset\reals^{\naturals}$, $E$ is countable dimensional if, and only if, $f\lbrack E\rbrack$ is first category.
\end{corollary}
{\bf Proof:} In Theorem \ref{erdsierp} take $\mathcal{I}$ to be {\sf CD}, the ideal of countable dimensional subsets of $\reals^{\naturals}$, and $\mathcal{J}$ to be $\mathcal{M}$, the ideal of first category subsets of $\reals^{\naturals}$.
By Hurewicz's theorem that $\reals^{\naturals}$ is not countable dimensional, ${\sf CD}$ is a $\sigma$-complete ideal on $\reals^{\naturals}$. By the Baire category theorem $\mathcal{M}$ is a $\sigma$-complete ideal on $\reals^{\naturals}$. 

Since $\reals^{\naturals}$ is separable, let $D$ be a countable dense set. Then $D$ is zero-dimensional. By Theorem \ref{tumarkin} there is a zero-dimensional ${\sf G}_{\delta}$ set $A\supset D$. But then $B=\reals^{\naturals}\setminus A$ is a first category set, and so (1) of Theorem \ref{erdsierp} is satisfied.  

Theorem \ref{tumarkin} also implies that each countable dimensional subset of $\reals^{\naturals}$ is contained in a countable union of zero-dimensional ${\sf G}_{\delta}$ sets, and so {\sf CD} has a cofinal family of cardinality $2^{\aleph_0}$. Since $\mathcal{M}$ has a cofinal family of cardinality $2^{\aleph_0}$, C.H. implies that (2) of Theorem \ref{erdsierp} holds. 
Evidently (3) and (4) also hold.
$\square$

\section{Hurewicz sets.}

Some strong analogies between the theory of the ideal of countable dimensional subsets of $\reals^{\naturals}$ and the theories of the ideal $\mathcal{N}$ of Lebesgue measure zero subsets of the real line, and the ideal $\mathcal{M}$ of first category subsets of uncountable complete separable metric spaces emerge when we consider the following features: Call a subset of $\reals^{\naturals}$ a \emph{Hurewicz set} if it is uncountable, but its intersection with every zero-dimensional subset of $\reals^{\naturals}$ is countable. Hurewicz sets are first category: Take a countable dense set $D\subset\reals^{\naturals}$. Since $D$ is zerodimensional it is, by Theorem \ref{tumarkin}, contained in a zerodimensional dense {\sf G}$_{\delta}$ set $E$. Since a Hurewicz set meets $E$ in only countably many points the Hurewicz set is first category. We shall see below that Hurewicz sets are first category in a strong sense.

The notion of a Hurewicz set is analogous to the notions of a Lusin set and of a Sierpinski set: A subset of $\reals$ is a \emph{Lusin set} if it is uncountable but its intersection with any first category subset of $\reals$ is countable. A subset of the real line is a \emph{Sierpinski set} if it is uncountable but its intersection with each Lebesgue measure zero set is countable. 

Mahlo \cite{mahlo} and Lusin \cite{lusin} independently showed that the Continuum Hypothesis implies the existence of a Lusin set of cardinality $2^{\aleph_0}$, and Sierpinski \cite{sierpinski24} showed that the Continuum Hypothesis implies the existence of a Sierpinski set of cardinality $2^{\aleph_0}$. It is known that the existence of a Lusin set of cardinality $2^{\aleph_0}$ does not imply the Continuum Hypothesis, and that the existence of a Sierpi\'nski set of cardinality $2^{\aleph_0}$ does not imply the Continuum Hypothesis. But Rothberger \cite{rothberger382} proved that the simultaneous existence of both a Lusin set of cardinality $2^{\aleph_0}$ and a Sierpinski set of cardinality $2^{\aleph_0}$ is equivalent to $2^{\aleph_0} = \aleph_1$. W. Hurewicz \cite{hurewicz32} showed that the existence of a Hurewicz set is equivalent to the Continuum Hypothesis. 

Some deeper properties of the Lusin sets and the Sierpi\'nski sets are tied up with the algebraic properties of the real line viewed as a topological group. 
\begin{theorem}[Galvin-Mycielski-Solovay]\label{gmslusin} If $L\subset \reals$ is a Lusin set, then for each first category set $M\subset \reals$, $L+M\neq \reals$.
\end{theorem}
Galvin, Mycielski and Solovay \cite{GMS} proved a significantly stronger theorem: Strong measure zero sets are characterized as sets whose translates by first category sets do not cover the real line. Theorem \ref{gmslusin} follows from this since Sierpi\'nski proved that Lusin sets are strong measure zero. Pawlikowski \cite{pawlikowski} proved the counterpart for Sierpi\'nski sets:
\begin{theorem}[Pawlikowski]\label{jpsierpinski} If $S\subset \reals$ is a Sierpi\'nski set, then for each Lebesgue measure zero set $N\subset\reals$, $S+N\neq\reals$.
\end{theorem}

We expect that a similar statement is true about Hurewicz sets in the topological group $(\reals^{\naturals},+)$:
\begin{conjecture}\label{conjecture}
If $H\subset \reals^{\naturals}$ is a Hurewicz set, then for each countable dimensional set $C\subset\reals^{\naturals}$, $H+C\neq \reals^{\naturals}$.
\end{conjecture} 

We have a partial result on this conjecture: Call a subset of $\reals^{\naturals}$ \emph{strongly countable dimensional} if it is a union of countably many closed, finite dimensional subsets. This notion was introduced by Nagata \cite{nagata} and Smirnov \cite{smirnov}. Every strongly countable dimensional set is countable dimensional, but not conversely. Thus, the $\sigma$-ideal {\sf SCD} generated by the strongly countable dimensional subsets of $\reals^{\naturals}$ is a proper subideal of {\sf CD}.
Indeed, every strongly countable dimensional set is of first category in $\reals^{\naturals}$, and thus {\sf SCD} is also a subideal of the $\sigma$-ideal of first category subsets of $\reals^{\naturals}$. Call a subset of $\reals^{\naturals}$ \emph{compactly countable dimensional} if it is a union of countably many compact finite dimensional sets. The $\sigma$-deal ${\sf KCD}$ generated by the compact finite dimensional sets is a proper subideal of ${\sf SCD}$.
We shall show
\begin{theorem}\label{scdsums} If $H$ is a Hurewicz subset of $\reals^{\naturals}$, then for every compactly 
countable dimensional set $C\subset\reals^{\naturals}$, $H+C$ is a first category subset of $\reals^{\naturals}$.
\end{theorem}

Since the union of countably many first category sets is first category, Theorem \ref{scdsums} follows directly from Lemma \ref{closedndim} below. This Lemma uses another important classical result:
\begin{theorem}[Hurewicz, Tumarkin]\label{hurewicztumarkin} In a separable metric space the union of countably many closed, $n$-dimensional sets, is an $n$-dimensional set.
\end{theorem}

\begin{lemma}\label{closedndim} If $H\subseteq \reals^{\naturals}$ is a Hurewicz set and $C\subset \reals^{\naturals}$ is a compact $n$-dimensional set, then $H+C$ is a first category subset of $\reals^{\naturals}$.
\end{lemma}
{\bf Proof:} Let $D\subset\reals^{\naturals}$ be a countable dense set. Then $D - C$ is a union of countably many closed $n$-dimensional sets, so by Theorem \ref{hurewicztumarkin} it is $n$-dimensional. By Theorem \ref{tumarkin} choose open sets $U_1\supset U_2\supset \cdots \supset U_k\supset \cdots\supset D - C$ such that $\bigcap_{k=1}^{\infty}U_k$ is $n$-dimensional. Then put $Y = \{x\in\reals^{\naturals}: x- C \subseteq \bigcap_{k=1}^{\infty}U_k\}$. Then $Y$ is dense in $\reals^{\naturals}$. Note that $Y$ is a ${\sf G}_{\delta}$-set: Consider a $k$, and an $x\in Y$. Since $C$ is compact there is an open neighborhood $V_k(x)$ of $x$ such that $V_k(x) - C\subset U_k$. But then $V_k=\cup_{x\in Y}V_k(x)$ is an open set containing $Y$, and $V_k-C\subset U_k$. It follows that $Y = \cap_{k<\infty}V_k$.

Next, put
\[
  X = H \cap(\bigcap_{k=1}^{\infty}U_k).
\]
Since $H$ is a Hurewicz set and $\bigcap_{k=1}^{\infty}U_k$ is $n$-dimensional, $X$ is a countable set. Define $Z = Y\setminus (X+C)$. Since $X+C$ is a union of countably many closed, $n$-dimensional sets, it is an $n$-dimensional ${\sf F}_{\sigma}$-set and thus first-category in $\reals^{\naturals}$. Consequently $Z$ is co-meager. But
\[
  Z = \{x\in Y:(x-C)\cap H = \emptyset\}\subseteq Y\setminus(H+C).
\]
Thus, $H+C$ is meager. $\square$

One can show that if there are Hurewicz sets, then there are ones satisfying the property in Conjecture \ref{conjecture}. There are several well-studied classes of weakly infinite dimensional subsets of $\reals^{\naturals}$ which are not countable dimensional.  In Example 1 we give, under CH, an example of a set $S$ in a very restrictive class of weakly infinite dimensional spaces (but not countable dimensional), and a Hurewicz set $H$, such that $S+H = \reals^{\naturals}$. The following is an unresolved weaker instance of Conjecture 1:
\begin{problem}\label{weakerversion} Is it true that if $C$ is a strongly countable dimensional set in $\reals^{\naturals}$ and $H$ is a Hurewicz set, then $H+C\neq \reals^{\naturals}$?
\end{problem}

\section{New classes of weakly infinite dimensional spaces.}

Many examples in infinitary dimension theory belong to classes of weakly infinite dimensional spaces introduced in \cite{babinkostova1}. Our examples in connection with Conjecture \ref{conjecture} are also in these classes. These classes are defined by applying selection principles to specific types of open covers. We now introduce these covers and selection principles.

\begin{center}{\bf Classes of open covers}
\end{center}

 For a given space $X$ the symbol $\op$ denotes the collection of all open covers of $X$. The following special types of open covers are relevant to our discussion:
\begin{itemize}
  \item{$\opfd$: This is the collection of all open covers $\mathcal{U}$ such that there is for each finite dimensional $F\subset X$ a $U\in\mathcal{U}$ with $F\subseteq U$, but $X$ is not a member of $\mathcal{U}$. (Thus, we are assuming $X$ is not finite dimensional.)}
  \item{$\opcfd$: This is the collection of all open covers $\mathcal{U}$ such that there is for each \emph{closed} finite dimensional $F\subset X$ a $U\in\mathcal{U}$ with $F\subseteq U$, but $X$ is not a member of $\mathcal{U}$. (Thus, we are assuming $X$ is not finite dimensional.)}  
  \item{$\opkfd$: This is the collection of all open covers $\mathcal{U}$ such that there is for each \emph{compact} finite dimensional $F\subset X$ a $U\in\mathcal{U}$ with $F\subseteq U$, but $X$ is not a member of $\mathcal{U}$. (Thus, we are assuming $X$ is not compact and finite dimensional.)} 
  \item{$\opk$: This is the collection of all open covers $\mathcal{U}$ such that there is for each \emph{compact} $F\subset X$ a $U\in\mathcal{U}$ with $F\subseteq U$, but $X$ is not a member of $\mathcal{U}$. (Thus, we are assuming $X$ is not compact.)} 
  \item{$\om$: This is the collection of $\omega$ covers of $X$. An open cover $\mathcal{U}$ of $X$ is an $\omega$-cover if $X$ is not a member of $\mathcal{U}$, but for each finite subset $F$ of $X$ there is a $U\in\mathcal{U}$ with $F\subseteq U$.} 
  \item{$\ga$: This is the collection of $\gamma$ covers of $X$. An open cover $\mathcal{U}$ of $X$ is a $\gamma$-cover if it is infinite and every infinite subset of $\mathcal{U}$ is a cover of $X$.} 
  \item{$\op^{gp}$: This is the collection of groupable open covers: An open cover $\mathcal{U}$ of a space $X$ is \emph{groupable} if there is a partition $\mathcal{U} = \bigcup_{n\in\naturals}\mathcal{F}_n$ into finite sets $\mathcal{F}_n$ that are disjoint from each other, such that for each $x\in X$ there is an $N$ with $x\in \bigcup\mathcal{F}_n$ for all $n\ge N$.} 
\end{itemize}

We have the following inclusion relations among these classes of open covers:
\[
  \opfd \subset \opcfd \subset \opkfd \subset \Omega \subset \op;\hspace{0.5in} \opk\subset \opkfd;\hspace{0.5in} \Gamma\subset\Omega.
\]

It is also worth noting that
if $X$ is a separable metric space then each $\mathcal{U}\in \opkfd$ has a countable subset $\mathcal{V}\in\opkfd$. The same is true about $\opk$.

\begin{center}{\bf Three selection principles.}
\end{center}
Let $\mathcal{A}$ and $\mathcal{B}$ be families of sets. The selection principle, $\sone(\mathcal{A},\mathcal{B})$ states: 
\begin{quote} For each sequence $(A_n:n\in\naturals)$ of elements of $\mathcal{A}$ there is a sequence $(B_n:n\in\naturals)$ such that for all $n$ we have $B_n\in A_n$, and $\{B_n:n\in\naturals\}\in\mathcal{B}$.
\end{quote}

The selection principle $\sfin(\mathcal{A},\mathcal{B})$ states:
\begin{quote} There is for each sequence $(A_n:n\in\naturals)$ of members of $\mathcal{A}$ a sequence $(B_n:n\in\naturals)$ of finite sets such that for each $n$ we have $B_n\subset A_n$ and $\bigcup_{n\in\naturals}B_n\in\mathcal{B}$.
\end{quote}

The selection principle $\Sc(\mathcal{A},\mathcal{B})$ states: 
\begin{quote} For each sequence $(A_n:n\in\naturals)$ of elements of $\mathcal{A}$, there is a sequence $(B_n:n\in\naturals)$ where each $B_n$ is a refinement of the collection of sets $A_n$, each $B_n$ is a pairwise disjoint family, and $\bigcup_{n\in\naturals}B_n$ is a member of $\mathcal{B}$.
\end{quote} 

In our context the families $\mathcal{A}$ and $\mathcal{B}$ will be classes of open covers of a topological space. The selection principle $\sone(\mathcal{A},\mathcal{B})$ has the following monotonicity properties: If $\mathcal{A}\subset \mathcal{C}$ and $\mathcal{B}\subset \mathcal{D}$, then we have the following implications:
\[
  \sone(\mathcal{C},\mathcal{B})\Rightarrow \sone(\mathcal{A},\mathcal{B}); \hspace{0.5in} \sone(\mathcal{A},\mathcal{B})\Rightarrow \sone(\mathcal{A},\mathcal{D}).
\]
Replacing $\sone$ with $\sfin$ or $\Sc$ in these implications give corresponding facts for the other two selection principles.
 
Special instances of these three selection principles appear in classical literature. For example: $\sone(\op,\op)$ is known as Rothberger's property and was introduced in the 1938 paper \cite{rothberger38}. $\sfin(\op,\op)$ is known as Menger's property and was introduced in the 1925 paper \cite{hurewicz25}. $\Sc(\op,\op)$ is known as property C, and was introduced in the 1978 paper \cite{AG}. $\Sc(\op,\op)$ is a selective version of Bing's property of screenability \cite{bing}, thus also known as \emph{selective screenability}. Selective screenability is an important property in infinitary dimension theory. Let $\op_2$ denote the collection of open covers $\mathcal{U}$ with $\vert\mathcal{U}\vert\le 2$. Then $\Sc(\op_2,\op)$ is equivalent to Alexandroff's notion of \emph{weakly infinite dimensional}.
A space which is not weakly infinite dimensional is said to be \emph{strongly infinite dimensional}.

\begin{center}{\bf The class $\sone(\opkfd,\op)$}\end{center}

Metrizable spaces with property $\sone(\opkfd,\op)$ seem important. First: Theorem 17 of \cite{babinkostova1} shows that for separable metric spaces $\sone(\opfd,\op)\Rightarrow \Sc(\op,\op)$. Monotonicity properties of the selection principles give: $\sone(\opkfd,\op)\Rightarrow\sone(\opcfd,\op)\Rightarrow\sone(\opfd,\op)$. Thus for metrizable spaces $\sone(\opkfd,\op)$ implies the important property $\Sc(\op,\op)$. Second: The following lemma shows that these spaces have the classical Menger property $\sfin(\op,\op)$:
\begin{lemma}\label{kofdop} $\sone(\opkfd,\op) \Rightarrow \sone(\opk,\op) \Rightarrow \sfin(\op,\op)$.
\end{lemma}
Third: Several central examples from the theory of infinite dimensional spaces are in this class. For example: R. Pol's example in \cite{rpol} which shows that a weakly infinite dimensional space need not be countable dimensional is in the class $\sone(\opkfd,\op)$. 

The preservation of infinitary dimension properties by products is not well understood yet. Whenever a new class of infinite dimensionality is defined it is of interest to know how this class behaves under various mathematical constructions, like products. We now make some remarks about products by spaces in the class $\sone(\opkfd,\op)$, using recent examples constructed by E. Pol and R. Pol in \cite{erpol}:

\begin{proposition}\label{opkfdproducts}
Assume Martin's Axiom.
\begin{enumerate}
  \item{Finite powers of sets in $\sone(\opkfd,\op)$ need not be in $\Sc(\op,\op)$.} 
  \item{The product of a space in $\sone(\opkfd,\op)$ with the space of irrational numbers need not be weakly infinite dimensional.} 
\end{enumerate}
\end{proposition}
To see this consider the following: In the proof of Proposition 4.1 of \cite{erpol} E. Pol and R. Pol construct, using Martin's Axiom, spaces $E_0$ and $E_1$ which have $\sone(\opkfd,\op)$. To see that these two spaces have this property we consider elements of the argument given in part (6) of the proof of \cite{erpol}, Proposition 4.1: For $i\in\{0,1\}$ the spaces $E_i$ are of the form $E\cup S_i$ where 
\begin{itemize}
  \item $S_i$ is the union of countably many compact, finite dimensional subspaces, and
  \item For each open set $W\supset S_i$ in $E_i$, $\vert E_i\setminus W\vert<2^{\aleph_0}$. 
\end{itemize}
Consider $E_i$. Let $(\mathcal{U}_n:n<\infty)$ be a sequence in $\opkfd$. Write $S_i$ as a union of compact finite dimensional subspaces $S^i_1\subset S^i_2\subset\cdots\subset S^i_n\subset\cdots$. For each $n$ choose a $U_{2\cdot n-1}\in\mathcal{U}_{2\cdot n-1}$ with $S^i_n\subseteq U_{2\cdot n-1}$. Then $W = \bigcup_{n<\infty}U_{2\cdot n-1}$ is an open set containing $S_i$, and so $E_i\setminus W$ has cardinality less than $2^{\aleph_0}$. But Martin's Axiom implies that a separable metric space of cardinality less than $2^{\aleph_0}$ has Rothberger's property $\sone(\op,\op)$ (Theorem 5 of \cite{FremlinMiller}). Thus, choose for each $n$ a $U_{2\cdot n}\in\mathcal{U}_{2\cdot n}$ such that $E_i\setminus W \subseteq\bigcup_{n<\infty}U_{2\cdot n}$. 

It follows that the space $X$ constructed there as the free union of $E_0$ and $E_1$ still has $\sone(\opkfd,\op)$. Recall that a metric space $(X,d)$ has the Haver property if there is for each sequence $(\epsilon_n:n<\infty)$ of positive real numbers a sequence $(\mathcal{U}_n:n<\infty)$ such that each $\mathcal{V}_n$ is a pairwise disjoint family of open sets, each of $d$-diameter less than $\epsilon_n$, and $\bigcup_{n<\infty}\mathcal{U}_n$ is a cover of $X$. Metrizable spaces with property $\Sc(\op,\op)$ have the Haver property in all equivalent metrics. Since \cite{erpol} shows that $X^2$ has the Menger property but under some metric $X^2$ does not have the Haver property, it follows that $X^2$ does not have property $\Sc(\op,\op)$. 

To see the second item: In Theorem 6.1 of \cite{erpol} it is pointed out that the product of the space $E_i$ with the subspace of irrational numbers of the real line is not even $\Sc(\op_2,\op)$, that is, is strongly infinite dimensional.

\begin{center}{\bf The class $\sone(\opkfd,\om)$} \end{center}

$\sone(\opk,\om)$ was considered in Section 2 of \cite{Panpav}. Theorem 2.1 of \cite{Panpav} implies:
\begin{lemma}[Pansera, Pavlovic]\label{panpavs1} A space has $\sone(\opk,\om)$ if, and only if, it has $\sone(\opk,\op)$ in all finite powers.
\end{lemma}

\begin{proposition}\label{kofdom} If a space has property $\sone(\opkfd,\om)$ then it has property $\sfin(\om,\om)$.
\end{proposition}
{\flushleft{\bf Proof:}} Let $X$ be a space with $\sone(\opkfd,\om)$. By monotonicity properties of $\sone(\cdot,\cdot)$ X then has $\sone(\opk,\om)$. Then by Lemma \ref{panpavs1} $X^n$ has $\sone(\opk,\op)$  for all finite $n$. By Corollary \ref{kofdop} $X^n$ has $\sfin(\op,\op)$ for all $n$. By Theorem 3.9 of \cite{coc2} this is equivalent to: $X$ has $\sfin(\om,\om)$.
$\square$

Using standard techniques one can prove: If a space has the property $\sone(\opkfd,\om)$, then it has the property $\sone(\opkfd,\op)$ in all finite powers. We don't know if the converse is true:
\begin{problem}\label{opkfdpowers} For metrizable space $X$ is it true that if each $X^n$ has property $\sone(\opkfd,\op)$ then $X$ has $\sone(\opkfd,\om)$?
\end{problem}

It is likely that the spaces $E_0$ and $E_1$ of \cite{erpol} mentioned above are in $\sone(\opkfd,\om)$, but we have not examined this. The space $X$ which is the free union of $E_0$ and $E_1$ is an example of a space in $\sone(\opkfd,\op)$ but not in $\sone(\opkfd,\om)$. Thus, under Martin's Axiom, $\sone(\opkfd,\om)$ is a proper subclass of $\sone(\opkfd,\op)$. Below we construct under CH an example (Example \ref{opkfdomproduct}) of a space which has the property $\sone(\opkfd,\om)$ and which is not countable dimensional. 

\begin{center}{\bf The class $\sone(\opkfd,\ga)$} \end{center}
$\sone(\opk,\Gamma)$ was explicitly considered in \cite{cdkm}, Section 3, and in \cite{dkm}, Section 3, where the second implication of Corollary \ref{opkfdgamma} is noted. 

\begin{corollary}\label{opkfdgamma} $\sone(\opkfd,\Gamma) \Rightarrow \sone(\opk,\Gamma) \Rightarrow \sfin(\om,\op^{gp})$.
\end{corollary}

Theorem 12 of \cite{cdkm} shows that each $\sigma$-compact space has the property $\sone(\opk,\Gamma)$. This is not true for $\sone(\opkfd,\Gamma)$: The Hilbert cube is not weakly infinite dimensional, but is a $\sigma$-compact space. 
Combining \cite{dkm} Theorem 10 and \cite{dkm} Proposition 13 gives the following result:
\begin{lemma}[Di Maio, Kocinac, Meccarriello]\label{opkgammapowers} If $X$ has $\sone(\opk,\ga)$, then for all $n$ also $X^n$ has this property.
\end{lemma}

The Hurewicz covering property is a strengthening of the Menger property $\sfin(\op,\op)$ and is defined as follows: Topological space $X$ has the \emph{Hurewicz covering property} if: For each sequence $(\mathcal{U}_n:n<\infty)$ of open covers of $X$ there is a sequence $(\mathcal{V}_n:n<\infty)$ of finite sets such that for each $n$, $\mathcal{V}_n\subseteq \mathcal{U}_n$, and for each $x\in X$, for all but finitely many $n$, $x\in\bigcup\mathcal{V}_n$. It was shown in \cite{coc7} that the Hurewicz covering proprty is equivalent to the selection principle $\sfin(\om,\op^{gp})$.
\begin{corollary}\label{powerssoneokfdgamma} If a space $X$ has property $\sone(\opkfd,\ga)$ then it has $\sfin(\om,\op^{gp})$ in all finite powers.
\end{corollary}

Using the technique in the proof of \cite{babinkostova1} Theorem 17 and using the techniques in \cite{babinkostova1} Lemma 6 one obtains the following two results: 
\begin{theorem} $\sone(\opfd,\op^{gp})\Rightarrow\Sc(\op,\op^{gp})$.
\end{theorem}

\begin{lemma}\label{powerskfdgamma} If a space $X$ has $\sone(\opkfd,\Gamma)$, then it has $\Sc(\op,\op)$ in all finite powers.
\end{lemma}
{\flushleft{\bf Proof:}} Let $X$ be a space with property  $\sone(\opkfd,\ga)$. By Lemma \ref{opkfdgamma} it has property $\sone(\opk,\ga)$, and by Lemma \ref{opkgammapowers} it has $\sone(\opk,\ga)$ in all finite powers, and thus the Hurewicz property in all finite powers. Then by Corollary 13 of \cite{babinkostova1}, all finite powers of $X$ has the property $\Sc(\op,\op)$.
$\square$

 Let $X$ be an infinite dimensional space which satisfies $\sone(\opkfd,\Gamma)$. By \cite{babinkostova1} and Corollary \ref{opkfdgamma}, all finite powers of $X$ have the Haver property, and indeed, the product of $X$ with any space having the Haver property again has the Haver property. It follows that $X\times \irrationals$ has the Haver property. We don't know if this product must be selectively screenable:

\begin{problem} If $X$ has $\sone(\opkfd,\ga)$ then does $X\times\irrationals$ have $\Sc(\op,\op)$?
\end{problem}

\section{Examples.}

We now describe two infinite dimensional examples that will be used to demonstrate a number of facts in connection with the preceding sections. In both cases we are interested in showing properties of Hurewicz sets, and thus we are confined to assuming the Continuum Hypothesis.\vspace{0.2in}

\begin{center}{\bf A ZFC+CH example of $\sone(\opkfd,\om)$.}
\end{center}
\vspace{0.2in}

Write: $\reals_n:=\{x\in\reals^{\naturals}:(\forall j>n)(x(n)=0)\}$ and $\reals_{\infty}:=\bigcup_{n=1}^{\infty}\reals_n$. Standard arguments prove Lemma \ref{rinfinityshifts} and Corollary \ref{rinfinitymanyshifts} below.

\begin{lemma}\label{rinfinityshifts} Assume $G\supset\reals_{\infty}$ is an open subset of $\reals^{\naturals}$ and $I_1,\, \cdots,\, I_n$ are closed intervals of positive length. For each $f\in\reals^{\naturals}$ there are an $m>n$ and closed intervals $I_{n+1},\, \cdots,\, I_m$ of positive length such that
\[
  I_1\, \times\, \cdots\, \times I_m\, \times\, \reals^{\naturals} \subseteq G-f.
\]
\end{lemma}

\begin{corollary}\label{rinfinitymanyshifts} If $(G_n:n<\infty)$ is a sequence of $G_{\delta}$ subset of $\reals^{\naturals}$ such that each $G_n$ contains $\reals_{\infty}$, and if $(f_n:n<\infty)$ is  sequence of elements of $\reals^{\naturals}$, then $\bigcap_{n=1}^{\infty} (G_n - f_n)$ contains a Hilbert cube.
\end{corollary}

\begin{example}[CH]\label{opkfdomproduct} There are subsets $X$ and $Y$ of $\reals^{\naturals}$ such that:
\begin{enumerate}
\item{$X$ and $Y$ satisfy $\sone(\opkfd,\om)$;}
\item{$X\setminus\reals_{\infty}$ and $Y\setminus\reals_{\infty}$ are Hurewicz sets;}
\item{$X$ and $Y$ have $\sone(\opkfd,\om)$ and $X\times Y$ not.}
\item{$X\cup Y$ is $\sone(\opkfd,\op)$, and not $\sone(\opkfd,\om)$.} 
\end{enumerate}
\end{example}

Fix the following enumerations:
\begin{itemize}
  \item$(f_{\alpha}:\alpha<2^{\aleph_0})$, the list of all elements of $\reals^{\naturals}$
  \item$(F_{\alpha}:\alpha<2^{\aleph_0})$, the list of all finite dimensional ${\sf G}_{\delta}$ subsets of $\reals^{\naturals}$
  \item$(L_{\alpha}:\alpha<2^{\aleph_0})$, the list of all ${\sf G}_{\delta}$ subsets of $\reals^{\naturals}$ containing $\reals_{\infty}$
  \item$((U^{\alpha}_n:n<\omega):\alpha<2^{\aleph_0})$ where for each $\alpha$ each element of $\reals_{\infty}$ is in $U^{\alpha}_n$ for all but finitely many $n$, and each $U^{\alpha}_n$ is open in $\reals^{\naturals}$.
\end{itemize}   
For each $\alpha<2^{\aleph_0}$ the set $T_{\alpha} = \bigcap_{n<\infty}(\bigcup_{m\ge n} U^{\alpha}_m)$
is a ${\sf G}_{\delta}$ set containing $\reals_{\infty}$.

Since $G_0 = T_0\cap L_0$ is a dense ${\sf G}_{\delta}$ set which contains $\reals_{\infty}$. Thus $G_0\cap (G_0-f_0)$ contains, by Corollay \ref{rinfinitymanyshifts}, a homeomorphic copy of the Hilbert cube. Choose 
\[
  x_0\in (G_0\cap (G_0-f_0)) \setminus ((\reals_{\infty}\cup F_0) \cup ((\reals_{\infty}\cup F_0)-f_0))
\]
and fix $y_0\in G_0$ with $x_0 = y_0-f_0$. Then define $S^1_0(0)=\{n: x_0\in U^0_n\}$ and $S^2_0(0)=\{n: y_0\in U^0_n\}$.

Assume we have $0<\alpha<\omega_1$, and that we have selected for each $\gamma<\alpha$ an $x_{\gamma}$ and a $y_{\gamma}$ and defined sets $S^1_{\nu}(\gamma)$ and $S^2_{\nu}(\gamma)$, $\nu\le\gamma$ such that:
\begin{enumerate}
  \item $\delta<\gamma<\alpha\Rightarrow S^1_{\gamma}(\delta) = S^2_{\gamma}(\delta)=\omega$; 
  \item $\gamma\le\nu<\delta<\alpha\Rightarrow$
        \begin{enumerate}  
          \item $S^1_{\gamma}(\delta)\subseteq^*S^1_{\nu}(\delta)$;  $S^2_{\gamma}(\delta)\subseteq^*S^2_{\nu}(\delta)$
          \item $S^1_{\gamma}(\delta),\, S^1_{\gamma}(\nu)\subseteq S^1_{\gamma}(0)$;  $S^2_{\gamma}(\delta),\, S^2_{\gamma}(\nu)\subseteq S^2_{\gamma}(0)$ are all infinite
        \end{enumerate}
  \item $\gamma\le \delta<\alpha\Rightarrow \{n: x_{\delta}\in U^{\gamma}_n\}\supseteq S^1_{\gamma}(\delta)$; $\{n: y_{\delta}\in U^{\gamma}_n\}\supseteq S^2_{\gamma}(\delta)$ 
  \item $\gamma\le \delta<\alpha\Rightarrow \{x_{\delta},\, y_{\delta}\}\subseteq G_{\gamma}\setminus(\reals_{\infty}\bigcup(\cup_{\nu\le \delta} F_{\nu})\bigcup\{x_{\nu}:\nu<\delta\}\bigcup\{y_{\nu}:\nu<\delta\})$;
  \item $\delta<\alpha\Rightarrow x_{\delta}= y_{\delta}-f_{\delta}$.
\end{enumerate}
Towards selecting $x_{\alpha}$ and $y_{\alpha}$ we consider two cases:

{\flushleft{\bf Case 1:} $\alpha=\beta+1$}\\
 For $\delta<\alpha$ define $S^1_{\alpha}(\delta) = \omega = S^2_{\alpha}(\delta)$. For $\gamma\le \beta$ we define
\[
  T^1_{\gamma} = \bigcap_{n=1}^{\infty}(\bigcup_{m\ge n,\, m\in S^1_{\gamma}(\beta)} U^{\gamma}_m),
\]
\[
  T^2_{\gamma} = \bigcap_{n=1}^{\infty}(\bigcup_{m\ge n,\, m\in S^2_{\gamma}(\beta)} U^{\gamma}_m),
\]
and for $\alpha$ we define
\[
  T^1_{\alpha} =  T^2_{\alpha} = \bigcap_{n=1}^{\infty}(\bigcup_{m\ge n} U^{\alpha}_m).
\]
and
\[
  G_{\alpha} = \bigcap_{\gamma\le\alpha}(L_{\gamma}\cap T^1_{\gamma}\cap T^2_{\gamma}).
\]
Then $G_{\alpha}$ is a ${\sf G}_{\delta}$ set containing $\reals_{\infty}$ and so by Corollary \ref{rinfinitymanyshifts} the set $G_{\alpha}\cap (G_{\alpha} - f_{\alpha})$ contains a Hilbert cube. But the set
\[
  B_{\alpha} = (\reals_{\infty}\bigcup(\cup_{\gamma\le\alpha}F_{\gamma})\bigcup\{x_{\gamma}:\gamma<\alpha\}\bigcup\{y_{\gamma}:\gamma<\alpha\}) - (\{f_{\gamma}:\gamma\le\alpha\}\cup\{\underline{0}\})
\]
is countable dimensional and thus does not contain a Hilbert cube. Choose $x_{\alpha}\in (G_{\alpha}\bigcap (G_{\alpha}-f_{\alpha}))\setminus B_{\alpha}$ and fix $y_{\alpha}\in G_{\alpha}$ with $x_{\alpha}=y_{\alpha}-f_{\alpha}$.

Then (4) holds at $\alpha$. Since for $\gamma<\alpha$ we have $x_{\alpha}\in T^1_{\gamma}$ and $y_{\alpha}\in T^2_{\gamma}$ the sets $S^1_{\gamma}(\alpha)=\{n\in S^1_{\gamma}(\beta): x_{\alpha}\in U^{\gamma}_n\}$ and $S^2_{\gamma}(\alpha)=\{n\in S^2_{\gamma}(\beta): y_{\alpha}\in U^{\gamma}_n\}$ are infinite. Similarly the sets $S^1_{\alpha}(\alpha) = \{n: x_{\alpha}\in U^{\alpha}_n\}$ and 
$S^2_{\alpha}(\alpha) = \{n: y_{\alpha}\in U^{\alpha}_n\}$ are infinite.

{\flushleft{\bf Case 2:} $\alpha$ is a limit ordinal.}\\
For $\gamma<\alpha$ choose infinite sets $F^1_{\gamma}$ and $F^2_{\gamma}$ so that for all $\beta<\alpha$, $F^1_{\gamma}\subseteq^* S^1_{\gamma}(\beta)$ and $F^2_{\gamma}\subseteq^* S^2_{\gamma}(\beta)$.

For $\gamma\le \beta$ we define
\[
  T^1_{\gamma} = \bigcap_{n=1}^{\infty}(\bigcup_{m\ge n,\, m\in F^1_{\gamma}} U^{\gamma}_m),
\]
\[
  T^2_{\gamma} = \bigcap_{n=1}^{\infty}(\bigcup_{m\ge n,\, m\in F^2_{\gamma}} U^{\gamma}_m),
\]
and for $\alpha$ we define
\[
  T^1_{\alpha} =  T^2_{\alpha} = \bigcap_{n=1}^{\infty}(\bigcup_{m\ge n} U^{\alpha}_m).
\]
and
\[
  G_{\alpha} = \bigcap_{\gamma\le\alpha}(L_{\gamma}\cap T^1_{\gamma}\cap T^2_{\gamma}).
\]
Since $G_{\alpha}$ is a ${\sf G}_{\delta}$ set containing $\reals_{\infty}$, Corollary \ref{rinfinitymanyshifts} implies that the set $G_{\alpha}\cap (G_{\alpha} - f_{\alpha})$ contains a Hilbert cube. The set
\[
  B_{\alpha} = (\reals_{\infty}\bigcup(\cup_{\gamma\le\alpha}F_{\gamma})\bigcup\{x_{\gamma}:\gamma<\alpha\}\bigcup\{y_{\gamma}:\gamma<\alpha\}) - (\{f_{\gamma}:\gamma\le\alpha\}\cup\{\underline{0}\})
\]
is countable dimensional and does not contain a Hilbert cube.

Choose $x_{\alpha}\in (G_{\alpha}\bigcap (G_{\alpha}-f_{\alpha}))\setminus B_{\alpha}$ and fix $y_{\alpha}\in G_{\alpha}$ with $x_{\alpha}=y_{\alpha}-f_{\alpha}$.

Then (4) holds at $\alpha$. Since for $\gamma<\alpha$ we have $x_{\alpha}\in T^1_{\gamma}$ and $y_{\alpha}\in T^2_{\gamma}$ the sets $S^1_{\gamma}(\alpha)=\{n\in F^1_{\gamma}: x_{\alpha}\in U^{\gamma}_n\}$ and $S^2_{\gamma}(\alpha)=\{n\in F^2_{\gamma}: y_{\alpha}\in U^{\gamma}_n\}$ are infinite. Similarly the sets $S^1_{\alpha}(\alpha) = \{n: x_{\alpha}\in U^{\alpha}_n\}$ and 
$S^2_{\alpha}(\alpha) = \{n: y_{\alpha}\in U^{\alpha}_n\}$ are infinite.

This completes the description of the process for choosing $x_{\alpha}$ and $y_{\alpha}$ for $\alpha<\omega_1$. Finally define
$X = \reals_{\infty}\bigcup\{x_{\alpha}:\alpha<\omega_1\}$ and $Y = \reals_{\infty}\bigcup\{y_{\alpha}:\alpha<\omega_1\}$.

{\flushleft{\bf Claim 1:} }\emph{If $\mathcal{U}$ is a collection of open subsets of $\reals^{\naturals}$ which contains an infinite subset $\mathcal{V}$, each infinite subset of which covers $\reals_{\infty}$, then there are countable sets $A_X$ and $A_Y$ such that $\mathcal{V}$ is an $\omega$ cover of $X\setminus A_X$ and of $Y\setminus A_Y$.}
{\flushleft{\bf Proof of Claim 1:}} For let $\mathcal{V}$ be such a subfamily of $\mathcal{U}$. We may assume that $\mathcal{V}$ is countable. Thus, we may assume that for some $\alpha$, fixed from now on, $\mathcal{V}$ is $(U^{\alpha}_n:n<\omega)$, as enumerated before the construction f $X$ and $Y$. Put $A_X=\{x_{\gamma}:\gamma\le \alpha\}$ and $A_Y=\{y_{\gamma}:\gamma\le \alpha\}$. 

Consider any finite subset $F\subset X\setminus A_X$. We may write $F=\{x_{\gamma_1},\cdots,x_{\gamma_n}\}$ where we have $\alpha<\gamma_1<\cdots <\gamma_n$. From the construction of $X$ and $Y$ we have $S^1_{\alpha}(\gamma_n)\subseteq^*\cdots\subseteq^* S^1_{\alpha}(\gamma_1)$. Choose $n$ large enough that
\[
  S^1_{\alpha}(\gamma_m)\setminus n \subseteq \cdots \subseteq S^1_{\alpha}(\gamma_1).
\]
Then for $k\in S^1_{\alpha}(\gamma_m)$ we have $F\subset U^{\alpha}_m$. A similar argument applies to $Y\setminus A_Y$.

{\flushleft{\bf Claim 2:} }\emph{$X$ and $Y$ both satisfy the selection principle $\sone(\opkfd,\om)$.}
{\flushleft{\bf Proof of Claim 2:}} We prove it for $X$. The proof for $Y$ is similar. Let $(\mathcal{U}_n:n<\infty)$ be a sequence of families of open sets of $\reals^{\naturals}$ such that each is an element of $\opkfd$ for $X$. Write $\naturals = \bigcup_{k<\infty} Z_k$ where each $Z_k$ is infinite and $Z_k\cap Z_m =\emptyset$ whenever we have $k\neq m$.

Starting with $Z_1$, choose for each $n\in Z_1$ a $U_n\in\mathcal{U}_n$ such that every infinite subset of $\{U_n:n\in Z_1\}$ covers $\reals_{\infty}$. Choose a countable set $X_1\subset X$ such that $\{U_n:n\in Z_1\}$ is an $\omega$-cover of $X\setminus X_1$.

Continuing with $Z_2$, choose for each $n\in Z_2$ a $U_n\in\mathcal{U}_n$ such that every infinite subset of $\{U_n:n\in Z_2\}$ covers $X_1\cup \reals_{\infty}$. Choose a countable set $X_2\subset X$ such that $\{U_n:n\in Z_2\}$ is an $\omega$-cover of $X\setminus X_2$. Note that we may assume $X_2\cap X_1 = \emptyset$.

With $X_1,\,\cdots,\, X_n$ selected countable subsets of $X$, and for each $j\le n$, $U_k\in\mathcal{U}_k$ selected for $k\in Z_j$ so that $(U_k:k\in Z_j)$ is an $\omega$-cover of $X\setminus X_j$, and every infinite subset of $(U_k:k\in Z_j)$ covers $X_1\cup\cdots\cup X_{j-1}\cup \reals_{\infty}$. Choose $U_k\in\mathcal{U}_k$ for $k\in Z_{n+1}$ so that each infinite subset of $\{U_k:k\in Z_{n+1}\}$ is a cover of $\reals_{\infty}\cup  X_1\cup\cdots X_n$ and by Claim 1 choose a countable $X_{n+1}\subset\reals^{\naturals}$ such that $\{U_k:k\in Z_{n+1}\}$ is an $\omega$-cover of $X\setminus X_{n+1}$.

Thus we obtain $U_k:k\in\naturals$ and countable subsets $X_n$ of $X$, disjoint from each other, such that for each $n$, $\{U_k:k\in\naturals\}$ is an $\omega$-cover of $X\setminus X_n$. Since for each finite subset $F$ of $X$ there is an $n$ with $F\cap X_n=\emptyset$, it follows that $\{U_k:k<\infty\}$ is an $\omega$-cover for $X$.

{\flushleft{\bf Claim 3:} }\emph{The sets $X\setminus \reals_{\infty}$ and $Y\setminus \reals_{\infty}$ both are Hurewicz sets.}
{\flushleft{\bf Proof of Claim 3:}} For consider any zero dimensional subset of $\reals^{\naturals}$. It is contained in a zero dimensional set of form $F_{\gamma}$, and $F_{\gamma}\cap (X\setminus\reals_{\infty}) \subset \{x_{\nu}:\nu<\gamma\}$, a countable set. 

{\flushleft{\bf Claim 4:} }\emph{$X\setminus \reals_{\infty} - Y\setminus \reals_{\infty} = \reals^{\naturals}$.}
{\flushleft{\bf Proof of Claim 4:}} This follows from (5) in the construction of $X$ and $Y$. 
{\flushleft{\bf Remarks:}}
{\flushleft{(1)}} $X$ and $Y$ \emph{are Menger in all finite powers:}\\
{\flushleft{\bf Proof:}} Claim 2 and Corollary \ref{kofdom}.
{\flushleft{(2)}} \emph{There is a Hurewicz set $H$ and an $\sone(\opkfd,\om)$ set $S$ with $H+S = \reals^{\naturals}$.}\\
{\flushleft{\bf Proof:}} Put $H = -(Y\setminus\reals_{\infty})$ and $S=X$. Apply Claims 2, 3 and 4.
{\flushleft{(3)}} $\sone(\opkfd,\om)$ \emph{is not preserved by finite products:}\\
{\flushleft{\bf Proof:}} 
Since the map $(x,y)\mapsto x-y$ is continuous, if $X\times Y$ is Menger, so is $X-Y = \reals^{\naturals}$. But the latter is not Menger. Thus $X\times Y$ is not Menger. By Corollary \ref{kofdop} $X\times Y$ is not $\sone(\opkfd,\op).$ 

{\flushleft{(4)}} $\sone(\opkfd,\om)$ \emph{is not preserved by finite unions:}\\
{\flushleft{\bf Proof:}} For the same reason, $X\cup Y$ is not Menger in all finite powers, even though $X\cup Y$ is Menger. It follows that $X\cup Y$ has the property $\sone(\opkfd,\op)$, but not the property $\sone(\opkfd,\om)$, and also not the property $\sone(\opk,\om)$.\vspace{0.2in}

\begin{center}{\bf A ZFC+CH example of $\sone(\opkfd,\Gamma)$}\end{center}
\vspace{0.2in}

The point with the following example is to demonstrate that spaces in this class need not be countable dimensional. We do not know if non-countable dimensional examples of metric spaces in this class can be obtained without resorting to hypotheses beyond the standard axioms of Zermelo-Fraenkel set theory.

Put $J=\{0\}\cup\lbrack\frac{1}{2},1\rbrack$. Define $J_m =\{x\in J^{\naturals}:(\forall n > m)(x(n) = 0)\}$ and $J_{\infty} = \bigcup_{m\in\naturals}J_m$.
Since each $J_m$ is compact and $m$-dimensional, $J_{\infty}$ is compactly countable dimensional and is a dense subset of $J^{\naturals}$.

For $X$ any subset of $\naturals$, define: $J(X) = \{x\in J^{\naturals}: \{n:x(n)\neq 0\}\subseteq X\}$.
  When $X$ is infinite, define $X^* = \{Y\subset\naturals:Y \mbox{ is infinite and }Y\setminus X \mbox{ is finite}\}$, and
$J^*(X) = \bigcup\{J(Y):Y\in X^*\}$.

\begin{lemma}\label{makinggamma} Let an infinite set $X\subset \naturals$, a countable set $\mathcal{C}\subset J^{\naturals}$, as well as a family $\mathcal{J}$ of open subsets of $J^{\naturals}$ be given so that relative to $\mathcal{C}\cup J_{\infty}$, $\mathcal{J}\in\opkfd$. Then there are an infinite $Y\subset X$, and a sequence $D_1,\, D_2,\, \cdots,\, D_n,\, \cdots$ in $\mathcal{J}$ such that $J_{\infty}\cup\mathcal{C}\cup J^*(Y)\subseteq \bigcup_{m\in\naturals}(\bigcap_{m\ge n}D_m)$.
\end{lemma}
{\bf Proof:} Enumerate $\mathcal{C}$ bijectively as $(c_n:n<\infty)$.
Choose $k_0 = \min(X)$. Then $J(\{1,\,\cdots,\, k_0\}) = J^{k_0}\times\{0\}^{\naturals}\subset J_{\infty}$ is compact and finite dimensional. Since $\mathcal{J}$ is in $\opkfd$, choose $D_1\in\mathcal{J}$ with $J_0\cup \{c_0\}\cup J(\{1,\cdots,k_0\})\subset D_1$.

For each $x\in J(\{1,\cdots,k_0\})$ choose a basic open set $B_x$ of the form $I_1(x)\times\cdots\times I_{k_0}(x)\times \{0\}^{k(x)}\times J^{\naturals} \subset D_1$. Since $J(\{1,\cdots,k_0\}$ is compact, choose a finite set $\{x_1,\, \cdots,\, x_m\}\in J(\{1,\cdots,k_0\})$ for which
\[
  J(\{1,\cdots,k_0\})\subseteq \bigcup_{i=1}^m B_{x_i}\subseteq D_1.
\]
Choose $k_1\in X$ with $k_1>\max\{k_0+k(x_i):1\le i\le m\}$. Then we have
\[
  J^{k_0}\times\{0\}^{k_1-k_0} \times J^{\naturals} \subseteq D_1.
\]
This specifies $k_0$, $k_1$ and $D_1$.

Now $J(\{1,\cdots,k_1\}) = J^{k_1}\times\{0\}^{\naturals}\subset J_{\infty}$ is a compact and finite dimensional subset of $J_{\infty}$. Since $\mathcal{J}$ is in $\opkfd$ for $\mathcal{C}\cup J_{\infty}$, choose a $D_2\in\mathcal{J}$ with $J_1\cup \{c_0,\, c_1\}\cup J(\{1,\cdots,k_1\})\subseteq D_2$. 
For each $x\in J(\{1,\cdots,k_1\})$ choose a basic open set $B_x$ of the form $I_1(x)\times\cdots\times I_{k_1}(x)\times \{0\}^{k(x)}\times J^{\naturals} \subset D_2$. Since $J(\{1,\cdots,k_1\}$ is compact, choose a finite set $\{x_1,\, \cdots,\, x_m\}\in J(\{1,\cdots,k_1\})$ for which
\[
  J(\{1,\cdots,k_1\})\subseteq \bigcup_{i=1}^m B_{x_i}\subseteq D_2.
\]
Choose $k_2\in X$ with $k_2>\max\{k_1+k(x_i):1\le i\le m\}$. Then we have
\[
  J^{k_1}\times\{0\}^{k_2-k_1} \times J^{\naturals} \subseteq D_2.
\]

This specifies $k_0$, $k_1$ and $k_2$, as well as $D_1$ and $D_2$.

Continuing in this way we find sequences 
\begin{enumerate}
\item{$k_0<k_1<\cdots <k_n<\cdots$ in $X$ and}
\item{$D_1,\, D_2,\, \cdots, D_n,\cdots \in\mathcal{J}$}  
\end{enumerate}
such that for $n\ge 1$ we have:

$J_n\cup \{c_0,\,\cdots,\, c_n\}\cup J(\{1,\cdots,k_n\})\subseteq D_{n+1}$ and $J^{k_n}\times\{0\}^{k_{n+1}-k_n}\times J^{\naturals} \subseteq D_{n+1}$.

Put $Y=\{k_n:0\le n<\infty\}$, an infinite subset of $X$.

Consider any infinite set $Z\subset\naturals$ with $Z\setminus Y$ finite. Choose $k_n$ so large that 
\[
  Z\setminus\{1,\cdots,k_n\} \subseteq Y.
\]
Now $J(Z)\subseteq J^{k_m}\times\{0\}^{k_{m+1}-k_m} \times J^{\naturals} \subseteq D_m$ for all $m>n$, and so
\[
  J(Z)\subseteq \bigcap_{m>n}D_m.
\]
It follows that $J_{\infty}\cup \mathcal{C}\cup J^*(Y)\subseteq \bigcup_{n<\infty}(\bigcap_{m>n} D_m)$.
$\square$

For each $Z\in Y^*$ the space $J(Z)$ has a subspace homeomorphic to a Hilbert cube, and so $J(Z)$ is strongly infinite dimensional. Moreover, $J^*(Z)\subseteq J^*(Y)$. 

Now we construct the example:
\begin{example}[CH]\label{okfdgamma} 
There is a set $X\subset\lbrack 0,1\rbrack^{\naturals}$ which is not countable dimensional, but it has the property $\sone(\opkfd,\ga)$.
\end{example}
{\bf Proof:} Let $(H_{\alpha}:\alpha<\omega_1)$ enumerate the finite dimensional ${\sf G}_{\delta}$-subsets of $J^{\naturals}$. Also, let $(\mathcal{J}_{\alpha}:\alpha<\omega_1)$ enumerate all countable  families of open sets, and let $(f_{\alpha}:\alpha<\omega_1)$ enumerate $^{\naturals}\naturals$. 

Now using CH recursively choose $X_{\alpha}\subset \naturals$ infinite and $x_{\alpha}\in J(X_{\alpha})$ such that:
\begin{enumerate}
\item{$\alpha<\beta \Rightarrow X_{\beta}\subseteq^* X_{\alpha}$,}
\item{If $\mathcal{J}_{\alpha}$ is an $\opkfd$-cover of $J_{\infty}\bigcup \{x_{\beta}:\beta < \alpha\}$, then $X_{\alpha}\in\lbrack \naturals \rbrack^{\aleph_0}$ is chosen such that for each $\beta<\alpha$, $X_{\alpha}\subset^* X_{\beta}$ and for some sequence $(D_n:n\in\naturals)$ in $\mathcal{J}_{\alpha}$ we have $J_{\infty}\cup\{x_{\beta}:\beta<\alpha\}\cup J^*(X_{\alpha}) \subseteq\bigcup_{n\in\naturals}(\bigcap_{m\ge n}D_m)$, and}
\item{For each $\beta$, $x_{\beta}\not\in \cup_{\gamma\le \beta}H_{\gamma}$.}
\end{enumerate}

Here is how we accomplish (2): First, choose an infinite set $X\subset \naturals$ such that for all $\beta<\alpha$, $X\subset^*X_{\beta}$. Since ${\mathcal{J}}_{\alpha}$ is an $\opkfd$-cover for $J_{\infty}\cup\{x_{\beta}:\beta<\alpha\}$, choose by Lemma \ref{makinggamma} an $X_{\alpha}\in\lbrack X \rbrack^{\aleph_0}$ and a sequence $(D_n:n\in\naturals)$ from ${\mathcal{J}}_{\alpha}$ such that $J^*(X_{\alpha}) \subseteq \bigcup_{n\in\naturals}(\bigcap_{m\ge n} D_m))$. Then we have 
\[
 J_{\infty}\cup J^*(X_{\alpha}) \cup \{x_{\beta}:\beta< \alpha\} \subseteq \bigcup_{n\in\naturals}(\bigcap_{m\ge n} D_m).
\]
Then as $J(X_{\alpha})$ contains a homeomorphic copy of the Hilbert cube, choose $x_{\alpha}\in J(X_{\alpha})\setminus \bigcup_{\beta\le\alpha}H_{\beta}$.

Put $Y = \{x_{\alpha}:\alpha<\omega_1\}$ and put $X = J_{\infty}\bigcup Y$. 
{\flushleft{\bf Claim 1:}} $Y$ \emph{is a Hurewicz set}.\\
For let $F$ be a finite dimensional subset of the Hilbert cube. Choose $\alpha<\omega_1$ so that $F\subset H_{\alpha}$, Then, as $Y\cap H_{\alpha}\subset \{x_{\beta}:\beta<\alpha\}$, it follows that $Y\cap F$ is countable. Thus $Y$ has a countable intersection with countable dimensional sets. 

{\flushleft{\bf Claim 2:}} \emph{Each open cover of $X$ which is in $\opkfd$ contains a subset which is a $\gamma$-cover for} $X$.\\
For let $\mathcal{U}$ be in $\opkfd$ for $X$. Choose an $\alpha$ with $\mathcal{U} = \mathcal{J}_{\alpha}$. Now at stage $\alpha$ we chose $X_{\alpha}$ and a sequence $(D_n:n\in\naturals)$ from $\mathcal{J}_{\alpha}$ so that $J^*(X_{\alpha})\bigcup J_{\infty}\bigcup \{x_{\beta}:\beta< \alpha\} \subseteq \bigcup_{n\in\naturals}(\bigcap_{m\ge n}D_m)$. Since for all $\gamma>\alpha$ we have $X_{\gamma}\subset^* X_{\alpha}$, it follows that for all $\gamma>\alpha$ we have $J^*(X_{\gamma})\subseteq J^*(X_{\alpha})$. But then $X\subseteq \bigcup_{n\in\naturals}(\bigcap_{m\ge n}D_m)$. 

{\flushleft{\bf Claim 3:}} \emph{$X$ has $\sone(\opkfd,\ga)$}. \\
For let $(\mathcal{U}_n:n<\infty)$ be a  sequence from $\opkfd$ for $X$. Choose an infinite discrete subset $(d_n:n<\infty)$ from $X$. For each $n$, put $\mathcal{V}_n=\{U\setminus\{d_n\}:U=V_1\cap\cdots\cap V_n \mbox{ and for }i\le n, \, V_i\in\mathcal{U}_i\}\setminus\{\emptyset\}$. Then put $\mathcal{U} = \bigcup_{n<\infty}\mathcal{V}_n$. Still $\mathcal{U}$ is in $\opkfd$ for $X$. By Claim 2 choose $D_m\in\mathcal{U}$, $m<\infty$, such that $X\subseteq\bigcup_{m\in\naturals}\bigcap_{n\ge m}D_n$. Choose $n_1$ so large that $d_1\in D_{m}$ for all $m\ge n_1$. Then choose $n_2>n_1$ so large that for an $m_1<n_2$, $D_{n_1}\in \mathcal{V}_{m_1}$ and $\{d_1,\cdots,d_{\max(\{n_1,m_1\}}\}\subseteq D_k$ for all $k\ge n_2$. Choose $n_3>n_2$ so large that for an $m_2<n_3$ we have $D_{n_2}\in\mathcal{V}_{m_2}$, and for all $k\ge n_3$ we have $\{d_j:j\le \max(\{n_2,m_2\})\}\subset D_k$, and so on.

We have $m_1<m_2<\cdots<m_k$, and $X\subseteq \bigcup_{k<\infty}\bigcap_{j\ge k}D_{n_j}$, and each $D_{n_j}$ is of the form $\bigcap_{i\le m_j} V^j_i$. For $i\le m_1$ put $U_i = V^1_i$, and for $k\ge 1$ and for $m_k<i\le m_{k+1}$ put $U_i = V^{k+1}_i$. Then we have for each $i$ that $U_i\in\mathcal{U}_i$, and $\{U_i:i<\infty\}$ is a $\gamma$-cover of $X$.

By Claim 1 $X$ is not countable dimensional. 
$\square$

{\flushleft{\bf Remarks:}}\\ 

When we choose $X_{\alpha}$ we may assume the enumeration function of $X_{\alpha}$ dominates $f_{\alpha}$. This guarantees that the subspace $Y$ of $X$ does not have the Menger property. 

By Corollary \ref{powerssoneokfdgamma} $X$ has the Hurewicz covering property in all finite powers. By Corollary 13 of \cite{babinkostova1} $X$ has the property $\Sc(\op,\op)$ in all finite powers.

In Theorem 15 of \cite{cdkm} it was shown that a Tychonoff space $X$ has the property $\sone(\opk,\ga)$ if, and only if, the space $C(X)$ consisting of continuous real-valued functions on $X$ has the following property: For every sequence $(A_n:n<\infty)$ of subsets of $X$, each having the function $f\in \overline{A_n}$ in the compact-open topology on $C(X)$, there is a sequence $(f_n:n\in\naturals)$ in $C(X)$ such that for each $n$, $f_n\in A_n$, and the sequence $(f_n:n<\infty)$ converges to $f$ in the point-open topology on $C(X)$. Since $\sone(\opkfd,\ga)$ implies $\sone(\opk,\ga)$, it follows that for the space $X$ constructed 
in example \ref{okfdgamma}, the function space $C(X)$ has the corresponding property.

\section{Added in proof:}

Shortly after we submitted our paper Roman Pol informed us that a certain modification of our construction of Example 2 answers our Problem 3. With Pol's kind permission we include his remarks here:

With CH assumed and with the notation established in Example 2 do the following: In addition to listing the $H_{\alpha}$, $\mathcal{J}_{\alpha}$ and $f_{\alpha}$ in the beginning of the construction, also list all open sets in $J^{\naturals}$ as $(G_{\alpha}:\alpha<\omega_1)$. 

At stage $\alpha$, when picking the point $x_{\alpha}$ we proceed as follows : if $J(X_{\alpha}) \setminus (G_{\alpha}\cup(\cup_{\beta\le\alpha}H_{\beta}))$ is nonempty, $x_{\alpha}$ is chosen from this set. Otherwise, $x_{\alpha}$ is an arbitrary point in $J(X_{\alpha})\setminus (\cup_{\beta\le\alpha}H_{\beta})$ if this is a subset of $G_{\alpha}$.   

{\flushleft{\bf Claim 4:}} \emph{$Y = \{ x_{\alpha}: \alpha < \omega_1 \}$ is not in $S_c(\op_2,\op)$}. \\
To this end repeat an argument from pages 90, 91 of \cite{rpol2}:  For $x\in J^{\naturals}$ write $supp(x)=\{n:x(n)\neq 0\}$. For each i, let $U_{i,0}$ ( respectively, $U_{i,1}$) be an open set in $J^{\naturals}$ consisting of points $x$ such that for some $j$, $\vert supp(x)\cap \{ 1, 2, ..., j \}\vert=i$, and  $0 < x(j) < 5/6$  ( respectively, $x(j) > 4/6$ ). A key observation ( which justifies (18) on page 90 ) is the following . Let $X = \{ n_1, n_2, ...\}$ with $n_1 < n_2 < ...$.   Then the trace of $U_{i,0}$ ( respectively, $U_{i,1}$) on $J(X)$ consists of the points $x$ in $J(X)$ with  $0 < x(n_i) < 5/6$ ( respectively, $x(n_i) > 4/6$). 

Then for each $i$ set $\mathcal{U}_i = \{ U_{i,0}, U_{i,1}\}$. 
{\flushleft{\bf Subclaim 4.1:}} If $X$ is any infinite subset of $\naturals$, the traces of $\mathcal{U}_i$ on $J(X)$ form a sequence in $J(X)$ witnessing that $J(X)$ does not have $\Sc(\op_2,\op)$. \\
For consider the closed subspace $C$ of $J(X)$ defined as follows: For $n\in X$ put $S_n=\lbrack \frac{1}{2},1\rbrack$, and for $n\not\in X$ put $S_n=\{0\}$. Then $C = \prod_{n=1}^{\infty}S_n$ is homeomorphic to $\lbrack \frac{1}{2}, 1\rbrack^X$.

For each $i$, define: 
\[ A_{i0} = \{x\in J(X): \vert supp(x)\cap\{1,\cdots,j\}\vert = i \mbox{ and }x(j) = \frac{1}{2}\}\]
and
\[ A_{i1} = \{x\in J(X): \vert supp(x)\cap\{1,\cdots,j\}\vert = i \mbox{ and }x(j) = 1\}.\]
Then $A_{i0}$ and $A_{i1}$ are disjoint closed subsets of $J(X)$ and are essentially opposite faces of the Hilbert cube $\lbrack \frac{1}{2}, 1\rbrack^X$. Moreover, $A_{ij}\subset U_{ij}$ and $A_{ij}\cap\overline{U}_{i(1-j)}=\emptyset$ for $j\in\{0,1\}$. Then the argument used for the Hilbert cube can be adapted and applied to the sequence of opposite pairs of faces $A_{i0}$ and $A_{i1}$ to show that the sequence $(\mathcal{U}_n:n<\infty)$ does not have a sequence of disjoint refinements that cover $J(X)$. Observe that if $C$ is a countable dimensional subset of $J(X)$, then the same is true about $J(X)\setminus C$. 

Now, return to the set $Y$: Assume, aiming at a contradiction, that $Y$ is an $\Sc(\op_2,\op)$ - space. Let $\mathcal{V}_i$ be a disjoint open collection in $J^N$ such that $\mathcal{V}_i$ refines $\mathcal{U}_i$ and the union $G$ of all elements of all $\mathcal{V}_i$ contains $Y$. Then $G = G_{\alpha}$ for some $\alpha<\omega_1$ and since $x_{\alpha}$ is in $G_{\alpha}$, by the rules for choosing $x_{\alpha}$ we have that for some countable dimensional set $C$, $J(X_{\alpha})\setminus C$ is contained in $G_{\alpha}$.  Therefore, the traces of $V_i$ on $J(X_{\alpha})\setminus C$ provide a forbidden sequence of open disjoint refinements of the traces of $U_i$ on $J(X_{\alpha})\setminus C$, a contradiction.   

{\flushleft{\bf Claim 5:}} \emph{With $X = J_{\infty}  \cup Y$, the product $X  \times \irrationals$ is not in $\Sc(\op_2,\op)$}. \\
Indeed, let $u : J \longrightarrow \{0,1\}$ take 0 to 0 and the interval [$\frac{1}{2}$,1] to 1, and let $f : J^{\naturals} \longrightarrow \{0,1\}^{\naturals}$ be the product map. Then $u$ and $f$ are continuous maps, implying f (the graph of f) is a closed subset of $J^{\naturals}\times \{0,1\}^{\naturals}$. 

Since the set $\irrationals =\{x\in \{0,1\}^{\naturals}: \vert\{n:x(n)=1\}\vert =\aleph_0\}$ is homeomorphic to the irrationals, and $f^{-1}(\irrationals) = J^{\naturals} \setminus J_{\infty}\supseteq Y$, we have $f\cap (X\times \irrationals) = f\cap (Y\times \irrationals) =\{(y,f(y)):y\in Y\}$ is a closed subset of $X\times\irrationals$. The map $g:Y\longrightarrow f\cap(Y\times\irrationals)$ defined by $g(y) = (y,f(y))$ is one-to-one, continuous, and its inverse is the projection onto the first coordinate, which also is continuous, and therefore $g$ shows that $Y$ is homeomorphic to $f\cap(Y\times\irrationals)$. By Claim 4 the closed subset $f\cap (X\times \irrationals)$ of $X\times \irrationals$ is not $\Sc(\op_2,\op)$, and thus $X\times\irrationals$ is not $\Sc(\op_2,\op)$.

\section{Acknowledgements}

We would like to thank Elzbieta Pol and Roman Pol for communicating their results to us, for a careful reading of our paper, and for remarks that improved the quality of the paper.

\end{document}